\documentclass[12pt,hidelinks]{article}
\usepackage[utf8]{inputenc}
\usepackage{hyperref}
\usepackage{amsmath, amssymb, url, float, graphicx, float}
\usepackage{fullpage}
\usepackage[small,bf]{caption}
\usepackage{cite}
\usepackage{xcolor}

\newcommand{\tmin}{\theta^\mathrm{min}}
\newcommand{\tmax}{\theta^\mathrm{max}}

\newcommand{\ii}{\mathbf{i}}
\newcommand{\complex}{\mathbf{C}}

\renewcommand\phi\varphi
\newcommand{\eps}{\varepsilon}
\newcommand{\real}{\mathbf{Re}}
\newcommand{\imag}{\mathbf{Im}}

\newcommand{\ones}{\mathbf 1}
\newcommand{\reals}{\mathbf{R}}

\newcommand{\symm}{{\mbox{\bf S}}}  

\newcommand{\tr}{\mathop{\bf tr}}
\newcommand{\diag}{\mathop{\bf diag}}

\newcommand{\dom}{\mathop{\bf dom}} 

\newcommand{\cf}{{\it cf.}}
\newcommand{\eg}{{\it e.g.}}
\newcommand{\ie}{{\it i.e.}}

\newcommand{\BEAS}{\begin{eqnarray*}}
\newcommand{\EEAS}{\end{eqnarray*}}
\newcommand{\BEA}{\begin{eqnarray}}
\newcommand{\EEA}{\end{eqnarray}}
\newcommand{\BEQ}{\begin{equation}}
\newcommand{\EEQ}{\end{equation}}
\newcommand{\BIT}{\begin{itemize}}
\newcommand{\EIT}{\end{itemize}}

\title{Bounds on Efficiency Metrics in Photonics}
\author{Guillermo Angeris\thanks{These authors contributed equally.}\and Theo Diamandis\footnotemark[1] \and Jelena Vu\v ckovi\'c \and Stephen Boyd}
\date{March 2022}

\begin{document} 
\maketitle 

\begin{abstract}
In this paper, we present a method for computing bounds for a variety of
efficiency metrics in photonics, such as the focusing efficiency or the mode
purity. We focus on the special case where the objective function can be
written as the ratio of two quadratic functions of the field and show that
there exists a simple semidefinite programming relaxation for this problem. We
provide a numerical example of bounding the maximal mode conversion purity for
a device of given size. This paper is accompanied by an open source
Julia package for basic simulations and bounds.
\end{abstract}
\section*{Introduction}
Traditionally, photonic devices were designed by a scientist or engineer
(whom we will call a designer) for a specific application. This designer would 
piece together components from a library to create a device for the desired 
task. While effective in practice, this process is time consuming, possibly
irrelevant to the final application of the design itself, and may produce 
designs that are far from optimal.
In an alternative approach to constructing devices, a designer 
specifies what they want while forfeiting control 
of how the device is constructed to an optimization algorithm. 
This optimization algorithm then attempts to find a device which maximizes the designer-specified 
performance metric---a mathematical objective function that outputs a number representing how well 
the design matches the desired specifications.
In photonics, this approach is called ``inverse design.''

\paragraph{Inverse design.} Photonic inverse design
\cite{
	kaoMaximizingBandGaps2005,
	jiaoDemonstrationSystematicPhotonic2005,
	lalau-keralyAdjointShapeOptimization2013,
	doryInversedesignedDiamondPhotonics2019,
	jiangGlobalOptimizationDielectric2019,
	yangInversedesignedNonreciprocalPulse2020
}
has been extremely successful in
finding photonic chips designs with very good practical performance 
when compared to designs generated by traditional methods. Still, there 
is an outstanding question of whether there exist designs with much
better performance. For simple devices, such as spherical lenses, a 
designer can find the optimal design with basic algebra and ray optics. 
However, for more complicated devices, finding the optimal design with respect
to some performance metric is an open research problem. As a result, designs
are usually found using heuristic methods in
practice~\cite{moleskyInverseDesignNanophotonics2018, angerisHeuristicMethodsPerformance2021}.

\paragraph{Bounds.} Given a design generated using a heuristic, it is natural to wonder how much 
better one could have done.
To answer this question, we need to determine a design's suboptimality with respect 
to some performance metric.
Recently, there has been a large amount of work 
in this area, attempting to find bounds of this form for a variety of metrics,
including mode volume~\cite{zhaoMinimumDielectricResonatorMode2020},
free space concentration~\cite{shimMaximalFreeSpaceConcentration2020},
integral overlap~\cite{moleskyHierarchicalMeanFieldOperator2020, kuangComputationalBoundsLight2020},
among many others~\cite{
millerFundamentalLimitsExtinction2014,
angerisComputationalBoundsPhotonic2019,
moleskyBoundsAbsorptionThermal2019,
shimFundamentalLimitsNearField2019,
michonLimitsSurfaceenhancedRaman2019,
trivediBoundsScatteringAbsorptionless2020,
moleskyFundamentalLimitsRadiative2020,
moleskyGlobalOperatorBounds2020
}. Additionally, the focusing objectives shown in~\cite{schab22}, released
well after the preprint of this article, are included as a special case
of the formulation presented here.

\paragraph{This paper.} In this paper, we extend the current bound formulations to include
objective functions which can be expressed as the ratio of two
quadratic functions of the field. This type of objective
includes a number of efficiency metrics such as the focusing
efficiency, the mode purity, among many others. We show a numerical example
of these bounds and also provide a set of simple open source packages
that can be used to compute bounds for a number of inverse design problems
whose objectives can be phrased as quadratics or the ratio of quadratics.

\section{The problem of maximizing efficiency}\label{sec:problem}
In the general photonic design problem, a designer must design a device that
maximizes some objective function $f$ of the fields $z$ by choosing from a range of
possible permittivities $\theta$ of a device at each point in space. (For
example, this might mean that the designer is only able to choose some permittivity
between that of air or silicon at each point in the design domain.) 

We will assume that the fields $z$ must
satisfy the \emph{electromagnetic wave equation}, which can be written as
\begin{equation}\label{eq:physics}
Az + \diag(\theta)z = b,
\end{equation}
for some linear operator $A$ and excitation $b$. In general, we will work with
a discretization of the fields and permittivities, such that $\theta, z, b \in
\reals^n$ are represented as real-valued $n$-vectors, that $f : \reals^n \to
\reals$ is a function mapping $z$ to a real number, while $A \in
\reals^{n\times n}$ is a real $n\times n$ matrix. We have assumed that $A$,
$z$, and $b$ are real in this case, but the complex case can be reduced to the
real one by separating it into its real and imaginary parts. (We will see an
explicit example of how to do this later in this paper.) As a rough guideline,
we may view equation~\eqref{eq:physics} as the linear-algebraic generalization
of
\[
    -\underbrace{\nabla \times \nabla \times}_{A} \underbrace{E}_z + \underbrace{\omega^2\mu_0\eps E}_{\diag(\theta)z} = \underbrace{-\ii \omega \mu_0 J}_{b},
\]
where the linear operator $A$ corresponds to a discretization of
$-\nabla\times\nabla\times\cdot$, the design parameters $\theta$ correspond to a
(scaled) discretization of the permittivities $\eps$, the fields $z$, of course,
correspond to the field $E$, and the excitation $b$ corresponds to the current
$-\ii\omega \mu_0 J$.

Because the designer is only allowed to choose materials whose parameters range
within some interval, we will write $\tmin_i \le \theta \le \tmax_i$ for $i=1,
\dots, n$. Without loss of generality, we will assume that $\tmax = -\tmin
= \ones$ since~\eqref{eq:physics} can always be rescaled such
that this is true. (See, \eg,~\cite[\S2.2]{angerisHeuristicMethodsPerformance2021} for more details.) The general
optimization problem the designer wishes to solve is then:
\begin{equation}\label{eq:general-main}
\begin{aligned}
	& \text{maximize} && f(z)\\
	& \text{subject to} && Az + \diag(\theta)z = b\\
    &&& -\ones \le \theta \le \ones. 
\end{aligned} 
\end{equation} 
Here the variables are the fields $z \in \reals^n$ and the permittivities 
$\theta \in \reals^n$, while the problem data are the matrix 
$A \in \reals^{n\times n}$ and the excitation $b \in \reals^n$. Note that this 
problem, as stated, is NP-hard~\cite[\S2.3]{angerisHeuristicMethodsPerformance2021}, 
so finding its optimal value, which we will call $p^\star$, is likely to be
computationally infeasible except for very small problems.

\paragraph{Efficiency metrics.}
A common problem in photonic design (and, more generally, in physical design)
is the problem of maximizing an efficiency metric. We say an objective is an
\emph{efficiency metric} whenever, for any $z \in \dom f$, we have that
\begin{equation}\label{eq:eff-definition}
0 \le f(z) \le 1,
\end{equation}
or, in other words, that the objective value for a feasible field $z$ is always
a number between 0 and 1. (We may, of course, replace the upper bound of $1$
with any finite number, say $v$, but this is the same as defining a new
objective function $\tilde f = (1/v)f$ which
satisfies~\eqref{eq:eff-definition}.) Note that there are some cases in which
the function $f$ might be unbounded from above (or below) and are therefore not
`efficiency metrics' in the sense specified here. Even in these cases, the
relaxation method we present will hold, but it is not guaranteed to return
points that are `reasonable'; \ie, the relaxation might give bounds which are
trivial. (We note that, in practice, we still expect the results to be relatively tight,
even without these guarantees.)

\paragraph{Ratio of quadratics.}
In many important cases in photonic design, efficiency metrics can be written as the ratio of
two quadratics in $z$, \ie,
\begin{equation}\label{eq:quad-ratio}
f(z) = \frac{z^TPz + 2p^Tz + r}{z^TQz + 2q^Tz + s},
\end{equation}
where $P, Q \in \symm^n$ are two symmetric matrices, while $p, q \in \reals^n$ and $r, s \in \reals$, whenever
$z^TQz + 2q^Tz + s > 0$ and is $-\infty$ otherwise. Note that this function $f$ is, in general, nonconvex.
(We will see some examples of such objective functions soon.) In order for $f$
to be an efficiency metric~\eqref{eq:eff-definition}, the numerator and
denominator must satisfy
\[
0 \le z^TPz + 2p^Tz + r \le z^TQz + 2q^Tz + s,
\]
for all $z \in \reals^n$. By minimizing over $z$, this is true whenever
\begin{equation}\label{eq:mat-inequality}
0 \le \begin{bmatrix}
    P & p\\
    p^T & r
\end{bmatrix} \le \begin{bmatrix}
    Q & q\\
    q^T & s
\end{bmatrix},
\end{equation}
where the inequalities are semidefinite inequalities~\cite[\S2.4.1]{cvxbook}.
The inequalities of~\eqref{eq:mat-inequality} imply that $P$ and $Q$ satisfy $0 \le
P \le Q$, while $r$ and $s$ must satisfy $0 \le r \le s$.

\paragraph{Optimization problem.} The resulting optimization problem, when $f$ is the ratio of
two quadratics, is:
\begin{equation}\label{eq:main}
\begin{aligned}
	& \text{maximize} && \frac{z^TPz + 2p^Tz + r}{z^TQz + 2q^Tz + s}\\
	& \text{subject to} && Az + \diag(\theta)z = b\\
	&&& -\ones \le \theta \le \ones.
\end{aligned}
\end{equation}
The variables in this problem are the fields $z \in \reals^n$ and the design
parameters $\theta\in \reals^n$, while the data are the matrices $A \in
\reals^{n\times n}$ and $P, Q \in \symm^n$, the vectors $p, q \in \reals^n$,
and the scalars $r, s \in \reals$. From the previous discussion, if the
objective is an efficiency metric, then the optimal value of~\eqref{eq:main},
$p^\star$, will also satisfy $0 \le p^\star \le 1$. Finding an upper bound to
this optimal value $p^\star$ would then give us an upper bound on the maximal
efficiency of the best possible design. 

\subsection{Example efficiency metrics}\label{sec:metrics}

\paragraph{Normalized overlap.} One important special case of an efficiency metric is sometimes known as the
\emph{normalized overlap}. The normalized overlap is defined as
\[
f(z) = \frac{(c^Tz)^2}{\|z\|_2^2},
\]
where $c \in \reals^n$ is a normalized vector with $\|c\|_2^2 = 1$. This is a special case of~\eqref{eq:quad-ratio}
where $P = cc^T$, $Q=I$, and $p = q = 0$, while $r = s = 0$.

It is easy to verify that this is indeed
an efficiency metric since $f(z) \ge 0$ as it is the ratio of two nonnegative quantities, while
\[
f(z) = \frac{(c^Tz)^2}{\|z\|_2^2} \le \frac{\|c\|_2^2\|z\|_2^2}{\|z\|_2^2} = \|c\|_2^2 = 1,
\]
where the first inequality follows from Cauchy--Schwarz~\cite[\S3.4]{linalg}. 
Whenever $c$ is a mode of the system, this objective is sometimes called the
normalized mode overlap, or the mode purity, and can be interpreted as the
fraction of power that is coupled into the mode specified by $c$, compared to
the total fraction of power going to all possible output modes.

In the case we wish to measure the normalized overlap only over some region
specified by indices $S \subseteq \{1, \dots, n\}$, we can instead write
\[
f(z) = \frac{(c^TRz)^2}{\|Rz\|_2^2},
\]
where the matrix $R \in \reals^{n\times n}$ is a diagonal matrix with diagonal entries
\begin{equation}\label{eq:selector}
R_{ii} = \begin{cases}
    1 & i \in S\\
    0 & \text{otherwise},
\end{cases} 
\end{equation}
for $i=1, \dots, n$. The resulting objective can be written in as the special
case of~\eqref{eq:quad-ratio} where $P = Rcc^TR$ and $Q = R^2 = R$, while $q =
p = 0$ and $r = s = 0$, and is also easily shown to be an efficiency metric.

\paragraph{Focusing efficiency.} While there are many ways of defining the
focusing efficiency of a lens, one practical definition is as the ratio of the
sum of intensities over two regions, written
\[
f(z) = \frac{\|R'z\|_2^2}{\|Rz\|_2^2}.
\]
Here, the matrices $R, R' \in \reals^{n\times n}$ are defined as
\[
R_{ii} = \begin{cases}
    1 & i \in S\\
    0 & \text{otherwise},
\end{cases} \qquad
R_{ii}' = \begin{cases}
    1 & i \in S'\\
    0 & \text{otherwise},
\end{cases}
\]
where, $S' \subseteq S \subseteq \{1, \dots, n\}$ are sets of indices over
which we sum the square of the field. In this case, we call $S$ the focusing
plane and $S'$ the focusing region or focal spot, which is usually chosen to be
approximately the full width at half maximum (FWHM) of the intensity along $S$.

This metric is nonnegative as it is the ratio of two nonnegative functions, and
satisfies $f(z) \le 1$ as $S' \subseteq S$.
We can write this as the special case of~\eqref{eq:quad-ratio} where $P = R'$ and $Q = R$,
while $p = q = 0$ and $r = s = 0$.

\section{Homogenization and bounds}
In this section, we will show a transformation of problem~\eqref{eq:main} which results in a quadratic
objective with an additional quadratic constraint, by introducing a new variable. We will then show
how to construct basic bounds using procedures similar to those
of~\cite{kuangComputationalBoundsLight2020, moleskyHierarchicalMeanFieldOperator2020, angerisHeuristicMethodsPerformance2021} and show a few simple extensions.

\subsection{Homogenized problem}
The main difficulty of constructing bounds for~\eqref{eq:main} is that the
fractional objective is difficult to deal with. We will first give a
`heuristic' derivation and show that it is always an upper bound to the
original problem. We then show that the converse is true: this new problem is
equivalent to the original when $A + \diag(\theta)$ is invertible for all
$-\ones \le \theta \le \ones$.

The main idea behind this method is to dynamically scale the input excitation, $b$, by some
factor $\alpha \in \reals$, such that the denominator is always equal to 1. To do this, we 
replace equation~\eqref{eq:physics} with one where the input $b$ is scaled, to get
\[
Ay + \diag(\theta)y = \alpha b.
\]
Here $y$ is a new variable we will call the \emph{scaled field} as we can write
$y = \alpha z$. Plugging this into the objective, assuming that $z$ is
feasible, we find that
\[
f(z) = f(y/\alpha) = \frac{(1/\alpha)^2y^TPy +2(1/\alpha) p^Ty + r}{(1/\alpha)^2y^TQy + 2(1/\alpha) q^Ty + s} 
= \frac{y^TPy +2\alpha p^Ty + \alpha^2r}{y^TQy + 2\alpha q^Ty + \alpha^2s}.
\]
We will then constrain the denominator to equal 1, which results in the
\emph{homogenized} problem
\begin{equation}\label{eq:homogenized}
\begin{aligned}
	& \text{maximize} && y^TPy + 2\alpha p^Ty + \alpha^2r\\
	& \text{subject to} && y^TQy + 2\alpha q^Ty + \alpha^2 s = 1\\
	&&& Ay + \diag(\theta)y = \alpha b\\
	&&&  -\ones \le \theta \le \ones, \quad \alpha \ge 0.
\end{aligned}
\end{equation}
The variables in this problem are the scaled field $y \in \reals^n$ and the
scaling factor $\alpha \in \reals$, while the problem data are the same as that
of the original problem~\eqref{eq:main}.

\paragraph{Upper bound.} We will now show that this new homogenized
problem~\eqref{eq:homogenized} is an upper bound to the original problem. More
specifically we will show that every feasible field $z$ and design parameters
$\theta$ for~\eqref{eq:main} has a feasible scaled field $y$, scaling factor
$\alpha > 0$, using the same design parameters $\theta$, with the same
objective value.

First, note that $z$ is feasible for~\eqref{eq:main}, by definition, if $f(z) > -\infty$,
\ie, if $z$ satisfies
\[
z^TQz + 2q^Tz + s > 0, \qquad (A + \diag(\theta))z = b,
\]
for some $-\ones \le \theta \le \ones$. Based on this choice of $z$, we will set
\[
\alpha = \frac{1}{\sqrt{z^TQz + 2q^Tz + s}}, \qquad y = \alpha z,
\]
and show that this choice of $\alpha$ and $y$ satisfies the constraints 
of~\eqref{eq:homogenized} with the same objective
value. Plugging this value into the first constraint of~\eqref{eq:homogenized}, we see that
\[
y^TQy + 2\alpha q^Ty + \alpha^2s = \alpha^2(z^TQz + 2q^Tz + s) = 1,
\]
while the second constraint has
\[
(A + \diag(\theta))y = \alpha (A + \diag(\theta))z = \alpha b.
\]
Finally, the objective satisfies:
\[
y^TPy + 2\alpha p^Ty + \alpha^2r = \alpha^2(z^TPz + 2 p^Tz + r) = \frac{z^TPz + 2 p^Tz + r}{z^TQz + 2q^Tz + s} = f(z),
\]
so the objective value for $y$ and $\alpha$ for problem~\eqref{eq:homogenized} is the same
as $f(z)$, the objective value for~\eqref{eq:main} with field $z$.

\paragraph{Equivalence.} We will show that, in fact, problem~\eqref{eq:homogenized}
and problem~\eqref{eq:main} are equivalent in the special case where $A + \diag(\theta)$ is
invertible for any choice of $-\ones \le \theta \le \ones$. (We note that the problems
have the same optimal value even in the case where the physics equation is not always invertible,
but invertibility usually holds in practice.) We've shown that every feasible field $z$
and design parameters $\theta$ have a corresponding scaled fields $y$, scaling parameter
$\alpha$ (with the same design parameters $\theta$). We will now show the converse: every
scaled field $y$ with scaling parameter $\alpha$ that is feasible for~\eqref{eq:homogenized}
has some corresponding field $z$ for~\eqref{eq:main} with the same objective value.
We break this up into two cases, one in which $\alpha \ne 0$ and one in which $\alpha = 0$.

Given $\alpha \ne 0$ and any $y$ satisfying the constraints of~\eqref{eq:homogenized},
we set $z = y/\alpha$.
This field $z$ satisfies the physics constraint with the same design parameters $\theta$ as
\[
(A + \diag(\theta))z = \frac{1}{\alpha}\left(A + \diag(\theta)\right)y = \frac{1}{\alpha}(\alpha b) = b.
\]
On the other hand, the objective value for this choice of $z$ is
\[
f(z) = f(y/\alpha) = \frac{y^TPy + 2\alpha p^Ty + \alpha^2r}{y^TQy + 2\alpha q^Ty + \alpha^2 s} = y^TPy + 2\alpha p^Ty + \alpha^2r.
\]
So this $z$ is also feasible with design parameters $\theta$ and the same objective value.

On the other hand, we will show that $\alpha = 0$ is never feasible for~\eqref{eq:homogenized}
for any choice of $-\ones \le \theta \le \ones$.
If $\alpha = 0$, then $(A + \diag(\theta))y = \alpha b = 0$. Since $A + \diag(\theta)$ is 
invertible by assumption, then $y=0$. This implies that
\[
y^TQy + 2\alpha q^Ty + \alpha^2 s = 0 \ne 1.
\]
So, given any $\theta$ and $\alpha = 0$, 
there is no scaled field $y$ that is feasible for~\eqref{eq:homogenized}. This shows that
the problems are equivalent as any feasible point for one is feasible in the other, with
the same objective value.

\subsection{Semidefinite relaxation}

In general, problem~\eqref{eq:homogenized} is still
nonconvex and likely computationally difficult to solve. On the other hand, we
can give a convex relaxation of the problem, yielding a new problem whose
optimal value is guaranteed to be at least as large as that
of~\eqref{eq:homogenized} while also being computationally tractable.

\paragraph{Variable elimination.} 
As in~\cite{angerisConvexRestrictionsPhysical2021, angerisHeuristicMethodsPerformance2021}, we 
can eliminate the design variable $\theta$ from problem~\eqref{eq:homogenized}, giving the following
equivalent problem over only the scaled field $y$ and scaling parameter $\alpha$,
\begin{equation}\label{eq:elim-design}
\begin{aligned}
	& \text{maximize} && y^TPy + 2\alpha p^Ty + \alpha^2r\\
	& \text{subject to} && y^TQy + 2\alpha q^Ty + \alpha^2s = 1\\
	&&& (a_i^Ty - \alpha b_i)^2 \le y_i^2, \quad i=1, \dots, n.
\end{aligned}
\end{equation}
with variables $y \in \reals^n$ and $\alpha \in \reals$. Here, $a_i^T$ denotes the $i$th
row of the matrix $A$, and the problem data are otherwise identical to that 
of~\eqref{eq:homogenized}. Additionally, we note that this problem is equivalent
to~\eqref{eq:homogenized} by the same argument as that of~\cite{angerisHeuristicMethodsPerformance2021}
and therefore to~\eqref{eq:main}.

\paragraph{Rewriting and relaxation.} The new problem~\eqref{eq:elim-design} is a nonconvex quadratically
constrained quadratic program (QCQP). We can write~\eqref{eq:elim-design} in a slightly more 
compact form:
\begin{equation}\label{eq:qcqp}
\begin{aligned}
	& \text{maximize} && x^T\bar Px\\
	& \text{subject to} && x^T\bar Q x = 1\\
	&&& x^T\bar A_i x \le 0, \quad i=1, \dots, n
\end{aligned}
\end{equation}
Here, the variable is $x = (y, \alpha) \in \reals^{n+1}$, while the problem data are the
matrices:
\[
\bar P = \begin{bmatrix}
    P & p\\
    p^T & r
\end{bmatrix}, \quad \bar Q = \begin{bmatrix}
    Q & q\\
    q^T & s
\end{bmatrix}, \quad \bar A_i = \begin{bmatrix}
a_ia_i^T - e_ie_i^T & -b_i a_i\\
-b_i a_i^T & b_i^2
 \end{bmatrix}, \quad i=1, \dots, n.
\]
Using this rewritten problem, we can then form a semidefinite relaxation in the following way:
\begin{equation}\label{eq:sdp}
\begin{aligned}
	& \text{maximize} && \tr(\bar P X)\\
	& \text{subject to} && \tr(\bar Q X) = 1\\
	&&& \tr(\bar A_i X) \le 0, \quad i=1, \dots, n\\
	&&& X \ge 0,
\end{aligned}
\end{equation}
where we are maximizing over the variable $X \in \symm^n$. We will call $d^\star$ the optimal 
value of this problem. Problem~\eqref{eq:sdp} is a
relaxation of~\eqref{eq:qcqp} as any feasible point $x \in \reals^n$
for~\eqref{eq:qcqp} gives a feasible point
$X = xx^T \ge 0$ for~\eqref{eq:sdp}, since
\[
\tr(\bar Q X) = \tr(\bar Q xx^T) = x^T\bar Q x = 1,
\]
with the same objective value, $\tr(\bar P X) = x^T\bar P x$.
This implies that the optimal objective value of~\eqref{eq:main},
$p^\star$ is never larger than the optimal objective value of~\eqref{eq:sdp};
\ie, we always have $p^\star \le d^\star$.


\paragraph{Properties.} There are several interesting basic properties of the
relaxation of problem~\eqref{eq:sdp}. First, since $\bar P \ge 0$ by
assumption~\eqref{eq:mat-inequality}, then $d^\star \ge 0$ since we know that, for any feasible 
$X$,
\[
d^\star \ge \tr(\bar PX) \ge 0.
\]
Since we also know from~\eqref{eq:mat-inequality} that $\bar P \le \bar Q$,
then, for any optimal $X^\star \ge 0$, we have
\[
d^\star = \tr(\bar P X^\star) \le \tr(\bar Q X^\star) = 1.
\]
This implies that
\[
0 \le p^\star \le d^\star \le 1,
\]
so $d^\star$ can always be interpreted as a percentage upper bound of
$p^\star$, as expected. We note that, even if $\bar P \le \bar Q$ does not
hold, the resulting problem~\eqref{eq:sdp} still yields a bound on the optimal
objective value $p^\star$. The difference is that we lose the guarantees
derived here that the resulting dual bound $d^\star$ satisfies $d^\star \le 1$.
Additionally, given any $X \ge 0$ with $\tr(\bar A_i X) \le 0$ for $i=1, \dots,
n$, and $\tr(\bar Q X) > 0$, then
\[
X^0 = \frac{1}{\tr(\bar Q X)}X
\]
is a feasible point for problem~\eqref{eq:sdp}.

Since we know that $\bar P \le \bar Q$, then the equality constraint $\tr(\bar
QX) = 1$, in problem~\eqref{eq:sdp} can be relaxed to $\tr(\bar QX) \le 1$,
with the same optimal objective value. Additionally, if we find a solution
$X^\star$ whose rank is 1, then $X^\star = xx^T$ for some $x$ and therefore we
have that $x = (y, \alpha)$ is a solution to the homogenized
problem~\eqref{eq:homogenized}, which is easily turned into a solution of the
original problem~\eqref{eq:main} by setting $z = y/\alpha$ and $\theta =
(a_i^Tz - b_i)/z_i$ when $z_i \ne 0$ and $0$ otherwise.

\paragraph{Dual problem.} The matrices $\bar A_i$ for $i=1, \dots, n$, $\bar Q$, and $\bar P$ 
are sometimes chordally-sparse~\cite{vandenbergheChordalGraphsSemidefinite2015}. This
structure can often be exploited to more quickly solve for the
optimal value of~\eqref{eq:sdp} by considering the dual problem instead. Applying
semidefinite duality~\cite[\S5.9]{cvxbook} to problem~\eqref{eq:sdp} gives
\begin{equation}\label{eq:sdp-dual}
	\begin{aligned}
		& \text{minimize} && \lambda_{n+1}\\
		& \text{subject to} && \sum_{i=1}^n \lambda_i \bar A_i + \lambda_{n+1}\bar Q \ge \bar P\\
		&&& \lambda \ge 0,
	\end{aligned}
\end{equation}
where $\lambda \in \reals^{n+1}$ is our optimization variable. This problem can then be
passed to solvers such as \texttt{COSMO.jl}~\cite{garstkaCOSMOConicOperator2019}, which support chordal 
decompositions, for faster solution times.

\paragraph{Discussion.} The transformation
of variables used here is very similar to the transformation used in the reduction
of linear fractional programs to linear programs~\cite[\S4.3.2]{cvxbook},
and similar transformations have been used for computational physics bounds
in~\cite{zhaoMinimumDielectricResonatorMode2020} in the special case that $b = 0$
and $Q = e_ie_i^T$
(see, \eg,~\cite[\S3.2]{angerisHeuristicMethodsPerformance2021}).
This family of variable transformations has been known in the optimization literature
since the 1960s~\cite{charnesProgrammingLinearFractional1962} for a specific subset of
optimization problems known as `fractional programming,' which include problems with objective
functions of the form of~\eqref{eq:quad-ratio}. The variable
transformation used on problem~\eqref{eq:main} to get the
homogenized problem~\eqref{eq:homogenized}
is sometimes called the generalized Charnes--Cooper transformation~\cite{schaibleParameterfreeConvexEquivalent1974}.
We also note that
the same methodology presented here can be applied to the formulation
in~\cite{angerisComputationalBoundsPhotonic2019, zhaoMinimumDielectricResonatorMode2020}, which is the special case where $P$ and
$Q$ are diagonal with nonnegative entries.
	
\subsection{Extensions}\label{sec:extensions}
There are a few basic extensions for the bounds provided in~\eqref{eq:sdp}.

\paragraph{Boolean constraints.}
If we are allowed to choose only Boolean parameters, \ie,
if we  have $\theta_i \in \{\pm 1\}$, instead of $-1 \le \theta_i \le 1$ for
each $i=1, \dots, n$, we can write the bound as
\[
\begin{aligned}
	& \text{maximize} && \tr(\bar PX)\\
    & \text{subject to} && \tr(\bar Q X) = 1\\
    &&& \tr(\bar A_i X) = 0, \quad i=1, \dots, n\\
    &&& X \ge 0,
\end{aligned}
\]
which follows from~\cite[\S3.2]{angerisHeuristicMethodsPerformance2021}. All of
the same properties for~\eqref{eq:sdp} also hold for the optimal value of this
problem.


\paragraph{Rewriting the physics equation.}
In practice, it is sometimes the case that the physics equation~\eqref{eq:physics} is 
better expressed in the following form:
\begin{equation}\label{eq:greens_formalism}
z + G\diag(\theta')z = b',
\end{equation}
where $0 \le \theta' \le \ones$, $b' \in \reals^n$, and $G \in \reals^{n\times n}$.
This formulation is sometimes called the `Green's formalism' or `integral equation' in 
electromagnetism and is equivalent to that of~\eqref{eq:physics}, in that
every $(z, \theta)$ that satisfies the physics equation~\eqref{eq:physics}
has a $\theta'$ such that $(z, \theta')$ satisfies~\eqref{eq:greens_formalism}, and vice versa.
To see this in the case that $A$ is invertible, we
can map~\eqref{eq:physics} to~\eqref{eq:greens_formalism} by setting $G = (2A - I)^{-1}$, $b' = Gb$,
and $\theta' = (\theta + \ones)/2$.

Similar to~\cite{kuangComputationalBoundsLight2020, moleskyHierarchicalMeanFieldOperator2020, angerisHeuristicMethodsPerformance2021}, we will reduce~\eqref{eq:greens_formalism},
which depends on both the field $z$ and the design parameters $\theta'$, to an equation depending only on
the \emph{displacement field} $w = \diag(\theta') z$. To do this, we can write~\eqref{eq:greens_formalism}
in terms of $w$ and $z$
\[
z + Gw = b', \quad w = \diag(\theta') z.
\]
Multiplying both sides of the first equation elementwise by $w$ gives:
\[
w_iz_i + w_ig_i^Tw = w_ib'_i, \quad i=1, \dots, n,
\]
where $g_i^T$ is the $i$th row of $G$.
Finally, because $0 \le \theta' \le 1$, we get that $w_i^2 = \theta_iw_iz_i \le w_iz_i$, which means that
\begin{equation}\label{eq:alternative}
w_i^2 + w_ig_i^Tw \le w_ib'_i, \quad i=1, \dots, n.
\end{equation}
The converse---that there exists a field $z$ and design parameters $\theta'$ satisfying~\eqref{eq:greens_formalism}
and $w = \diag(\theta')z$, for any $w$ satisfying~\eqref{eq:alternative}---can be easily shown;
\cf, \cite[App.\ A]{angerisHeuristicMethodsPerformance2021}.

Rewriting~\eqref{eq:greens_formalism} we have that $z = b' - Gw$, and replacing the physics constraint in~\eqref{eq:main} with~\eqref{eq:alternative}
gives a new problem over the displacement field $w$,
\[
\begin{aligned}
	& \text{maximize} && \frac{w^TP'w + 2p'^Tw + r'}{w^TQ'w + 2q'^Tw + s'}\\
	& \text{subject to} && w_i^2 + w_ig_i^Tw \le w_ib'_i, \quad i=1, \dots, n,
\end{aligned}
\]
with variable $w \in \reals^n$ and problem data $G$, $b$, and
\[
P' = G^TPG, \quad p' = -G^TP(p+b), \quad r' = b^TPb + 2p^Tb + r,
\]
while
\[
Q' = G^TQG, \quad q' = -G^TQ(q+b), \quad s' = b^TQb + 2q^Tb + s.
\]
Applying the same homogenization procedure and semidefinite relaxation, this results in a problem identical to~\eqref{eq:sdp}
with the following problem data:
\[
\bar P = \begin{bmatrix}
    P' & p'\\
    p'^T & r'
\end{bmatrix}, \quad \bar Q = \begin{bmatrix}
    Q' & q'\\
    q'^T & s'
\end{bmatrix}, \quad \bar A_i = \begin{bmatrix}
e_ie_i^T + (e_ig_i^T + g_ie_i^T)/2 & -b_i' e_i\\
-b_i' e_i^T & 0
 \end{bmatrix}, \quad i=1, \dots, n.
\]

\paragraph{Convex constraints.}
We can also allow convex constraints in the SDP relaxation~\eqref{eq:sdp}. If we have a number of
convex constraints on the field $z=y/\alpha$ given by $f_j: \reals^n \to \reals$ for $j=1, \dots, m$, we can write
\[
\begin{aligned}
	& \text{maximize} && \tr(\bar PX)\\
    & \text{subject to} && \tr(\bar Q X) = 1\\
    &&& \tr(\bar A_i X) = 0, \quad i=1, \dots, n\\
    &&& \alpha f_j\left(\frac{y}{\alpha}\right) \le 0, \quad j=1, \dots, m\\
    &&& X = \begin{bmatrix}
        Y & y\\
        y^T & \alpha
    \end{bmatrix}\ge 0.
\end{aligned}
\]
The variables in this problem are the matrices $X \in \symm^{n+1}$, $Y \in \symm^{n}$, the vector $y\in \reals^n$, and
scalar $\alpha \in \reals$, while the problem data are identical to that of~\eqref{eq:sdp}. This new problem is again
a convex optimization problem since the functions $\alpha f_j(y/\alpha)$ over the variable $(y, \alpha)$ are convex if
the original functions $f_j$ are convex. This transformation is known as the perspective transform and always preserves convexity~\cite[\S3.2.6]{cvxbook}.
The resulting problem is then convex and can therefore be efficiently solved in most cases.

\paragraph{Additional quadratic constraints.} Similar to the previous, we can include additional (potentially indefinite)
quadratic constraints on the field $z$ into the relaxation~\eqref{eq:sdp}. More specifically, we wish to include
a number of constraints on the field $z$,
\[
z^TU_jz + 2u_j^Tz + t_j \le 0,
\]
with matrices $U_j \in \symm^n$, vectors $u_j \in \reals^n$, and scalars $t_j \in \reals$ for
$j=1, \dots, m$. Using the fact that $z = y/\alpha$, we can write these as
\[
y^TU_iy + 2\alpha u_j^Ty + \alpha^2t_j \le 0, \quad j=1, \dots, m,
\]
or, equivalently as
\[
x^T\bar U_jx \le 0, \quad i=1, \dots, m,
\]
where $x = (y, \alpha)$ as in~\eqref{eq:sdp} and
\[
\bar U_j = \begin{bmatrix}
    U_j & u_j\\
    u_j^T & t_j
\end{bmatrix}, \quad j=1, \dots, m.
\]
Using the same relaxation method as in~\eqref{eq:sdp} with the additional
quadratic inequalities,
we get the following semidefinite problem:
\[
\begin{aligned}
	& \text{maximize} && \tr(\bar PX)\\
    & \text{subject to} && \tr(\bar Q X) = 1\\
    &&& \tr(\bar A_i X) = 0, \quad i=1, \dots, n\\
    &&& \tr(\bar U_j X) \le 0, \quad j=1, \dots, m\\
    &&& X \ge 0.
\end{aligned}
\]
This problem has the same variables and problem data as~\eqref{eq:sdp}, with the addition of the
matrices $\bar U_j \in \symm^{n+1}$, as defined above.

\section{Numerical experiments}
In this section, we solve problem~\eqref{eq:sdp} for the 
the maximal mode purity of a small mode converter. We also find a design that approximately
saturates the bound.
To compute these bounds, we introduce two open source
Julia~\cite{bezansonJuliaFreshApproach2017, legatMathOptInterfaceDataStructure2021} packages, 
\texttt{WaveOperators.jl} and \texttt{PhysicalBounds.jl}, that
allow users to setup
physical design problems and compute bounds in only a few lines of code.

Our packages setup the dual form of the SDP~\eqref{eq:sdp-dual} using 
JuMP~\cite{dunningJuMPModelingLanguage2017, legatMathOptInterfaceDataStructure2021}
and solve it using any conic solver that supports semidefinite programming.
We use \texttt{SCS}~\cite{odonoghueConicOptimizationOperator2016} for the 
experiments in this paper. The code can be found at
\begin{center}
\texttt{github.com/cvxgrp/WaveOperators.jl}\\
\vspace{.5em}
    \texttt{github.com/cvxgrp/PhysicalBounds.jl}
\end{center}
which can be used to generate the plots found in this paper.

\subsection{General physics set up}

\paragraph{Physics equation.} We assume that the EM wave equation is appropriately discretized and results in a problem
of the form
\[
Az +\diag(\theta)z = b.
\]
Here $z \in \complex^n$ is the (complex) field while $\theta\in\reals^{n}$ are the (real) parameters and $A \in \complex^{n\times n}$, $b \in \complex^n$. To turn this into a problem
over real variables, we can separate the real and imaginary parts of the variables to get a new physics equation that is purely real:
\[
A'z' + \diag(\theta, \theta)z' = b'.
\]
Here, we define:
\[
A' = \begin{bmatrix}
    \real(A) & -\imag(A)\\
    \imag(A) & \real(A)
\end{bmatrix}, \qquad
b' = \begin{bmatrix}
    \real(b)\\
    \imag(b)
\end{bmatrix}, \qquad
z'= \begin{bmatrix}
    \real(z)\\
    \imag(z)
\end{bmatrix},
\]
where $\real(x)$ denotes the elementwise real part of $x$ (where $x$ is a vector or a matrix) while $\imag(x)$ denotes the imaginary part.
Note that this results in a larger system with parameters $A' \in \reals^{2n \times 2n}$, $b' \in \reals^{2n}$, and field $z' \in \reals^{2n}$,
whose parameters are all real. Finally, note that we can write this system as
\[
A'z' + \diag(\theta')z' = b', \qquad \theta'_{n+i} = \theta_i',
\]
where we have introduced a new, larger vector of parameters, $\theta' \in \reals^{2n}$ with an additional
constraint. Dropping this latter constraint over $\theta'$ leads to a relaxation of the original physics equation, in the following
sense: any design and field that satisfies the original equation also satisfies this new `relaxed' equation. This makes the final
physics equation:
\begin{equation}
    Az + \diag(\theta)z = b,
\end{equation}
where we have dropped the apostrophes for convenience. As a reminder we have
the physics operator $A \in \reals^{2n\times 2n}$, excitation $b \in
\reals^{2n}$, the field $z \in \reals^{2n}$, and the permittivities $\theta \in
\reals^{2n}$. This relaxation corresponds to allowing the designer to vary both
real and imaginary permittivities, where each component is box-constrained,
while the original problem only allows the designer to choose real
permittivities. (We note that some solvers, including
\texttt{Hypatia.jl}~\cite{coey2021solving}, support complex variables, but we
do not solve the problem over complex variables in this work.)

\begin{figure}
    \centering
	\includegraphics[width=\linewidth]{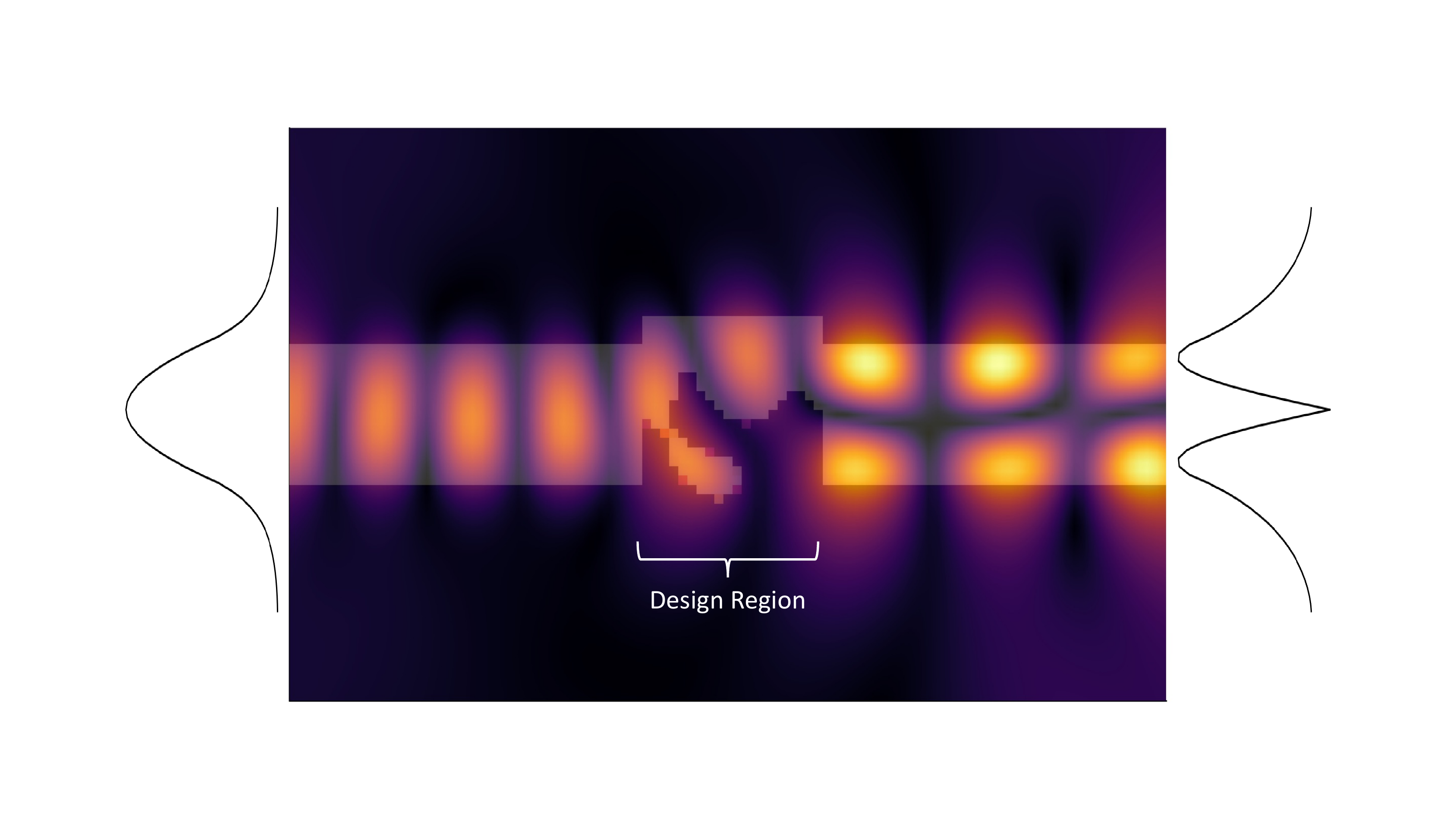}
	\caption{The designer wishes to choose materials in the design region to maximize the mode purity, measured at the 
    output of the waveguide.}
	\label{fig:setup-mode}
\end{figure}

\subsection{Mode converter}\label{sec:mode-example}
The setup is shown in figure~\ref{fig:setup-mode}. In this problem, the
designer is attempting to design a mode converter with the maximum mode purity,
by choosing the permittivities in the region shown. The input to this device is
the first order mode of the waveguide on the left hand side. The desired output
is a field whose normalized overlap with the second order mode of the waveguide
is maximized. In this problem, the designer is allowed to choose the
permittivities within the design region, so long as the permittivities lie in a
given interval. More information about the problem set up is given in
appendix~\ref{app:set-up} and the documentation of the corresponding packages.

\paragraph{Problem data.} In our specific problem set up,
as shown in figure~\ref{fig:setup-mode}
we have a source that is a distance of about one wavelength
from the design region. The simulation region is a rectangle that is one
wavelength tall and 1.6 wavelengths wide.
The design region is a centered square with side length 1/3 of a wavelength.
In this approximation, we assume that the grid is a $60\times 96$ grid; \ie,
the side length of a pixel in this simulation
is roughly $1/60$th of a free-space wavelength,
so $h = 1/60$. 
The material contrast (see appendix~\ref{app:set-up})
is set to $\delta = 10$
while the free-space wavenumber is $k = 2\pi$.

\paragraph{Optimization problem.}
In this experiment, we attempt to maximize the normalized overlap as defined
in~\S\ref{sec:metrics}:
\[
\begin{aligned}
	& \text{maximize} && \frac{(c^TRz)^2}{\|Rz\|_2^2}\\
	& \text{subject to} && Az + \diag(\theta)z = b\\
	&&& -\ones \le \theta \le \ones.
\end{aligned}
\]
Here the variables and problem data are similar to those of problem~\eqref{eq:main}.
More specifically, the problem variables are $z \in \reals^{2n}$, $\theta \in \reals^{2n}$, while the problem data
is the physics matrix $A \in \reals^{2n \times 2n}$, the excitation $b \in \reals^{2n}$,
the vector $c \in \reals^{2n}$ specifying the desired output mode,
and the matrix $R \in \reals^{2n \times 2n}$, defined in~\eqref{eq:selector}, where the region $S$ is the rightmost
column of pixels.
The resulting semidefinite upper bound for this problem
is given in~\eqref{eq:sdp} with
\[
P = Rcc^TR, \qquad Q = R, \qquad p=0, \qquad q=0, \qquad r = 0, \qquad s = 0.
\]

\begin{figure}[t]
    \centering
	\includegraphics[width=0.45\linewidth]{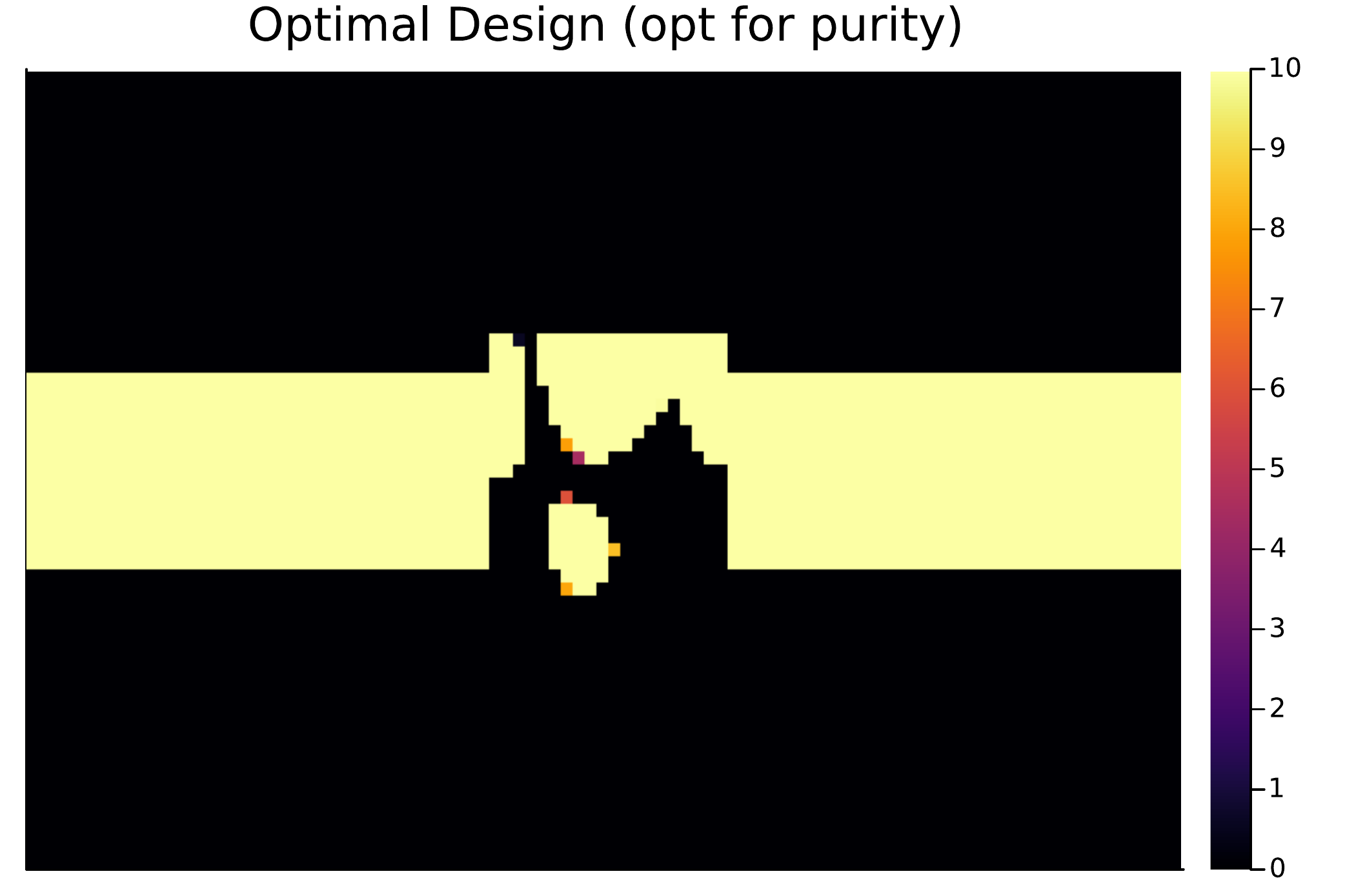}
    \includegraphics[width=0.45\linewidth]{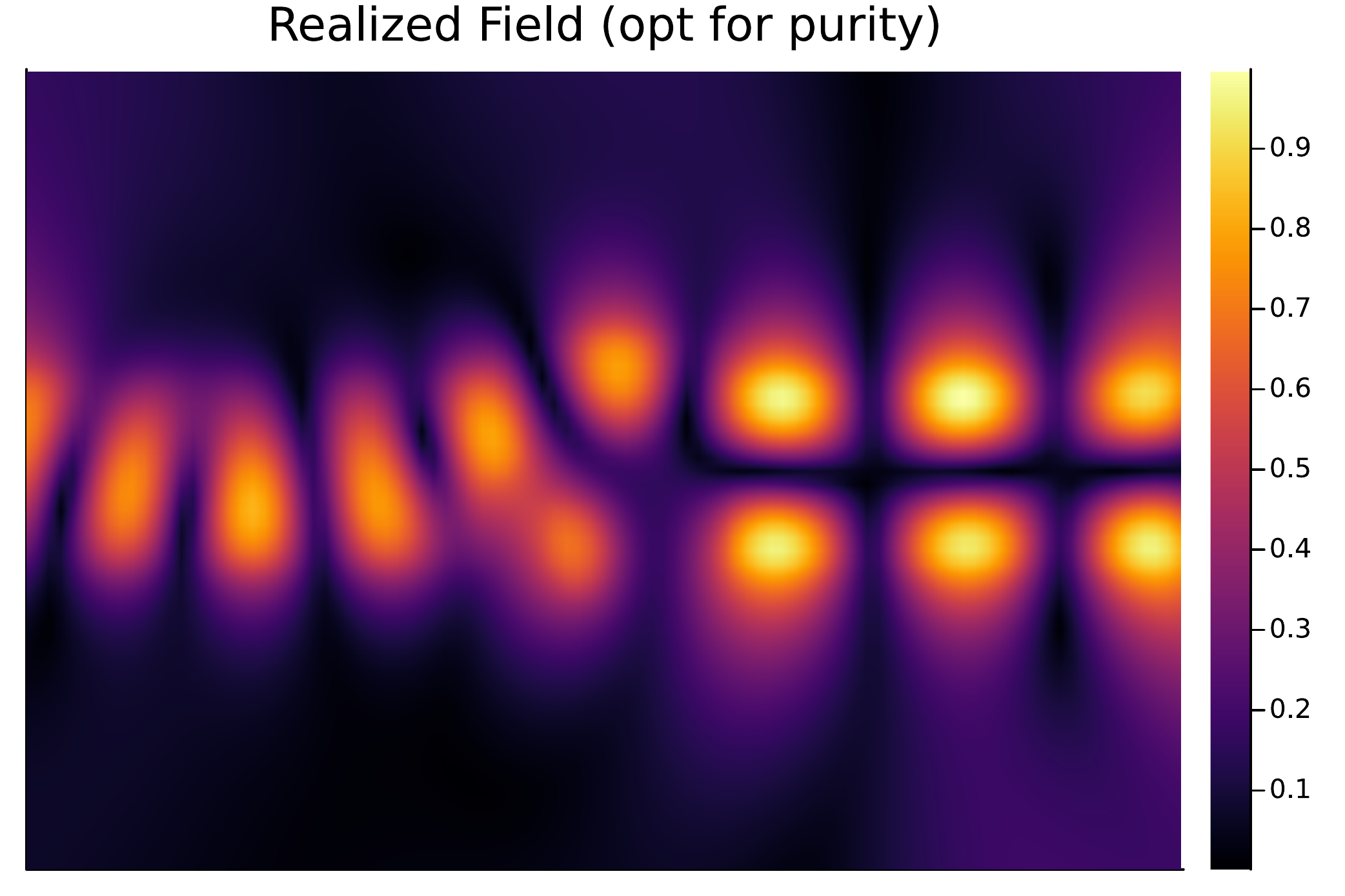}
	\caption{The design (left) is optimized for mode purity. The corresponding field 
    (right) closely matches the target mode at the output.}
	\label{fig:purity-deisgn-field}
\end{figure}

\begin{figure}[t]
    \centering
	\includegraphics[width=0.45\linewidth]{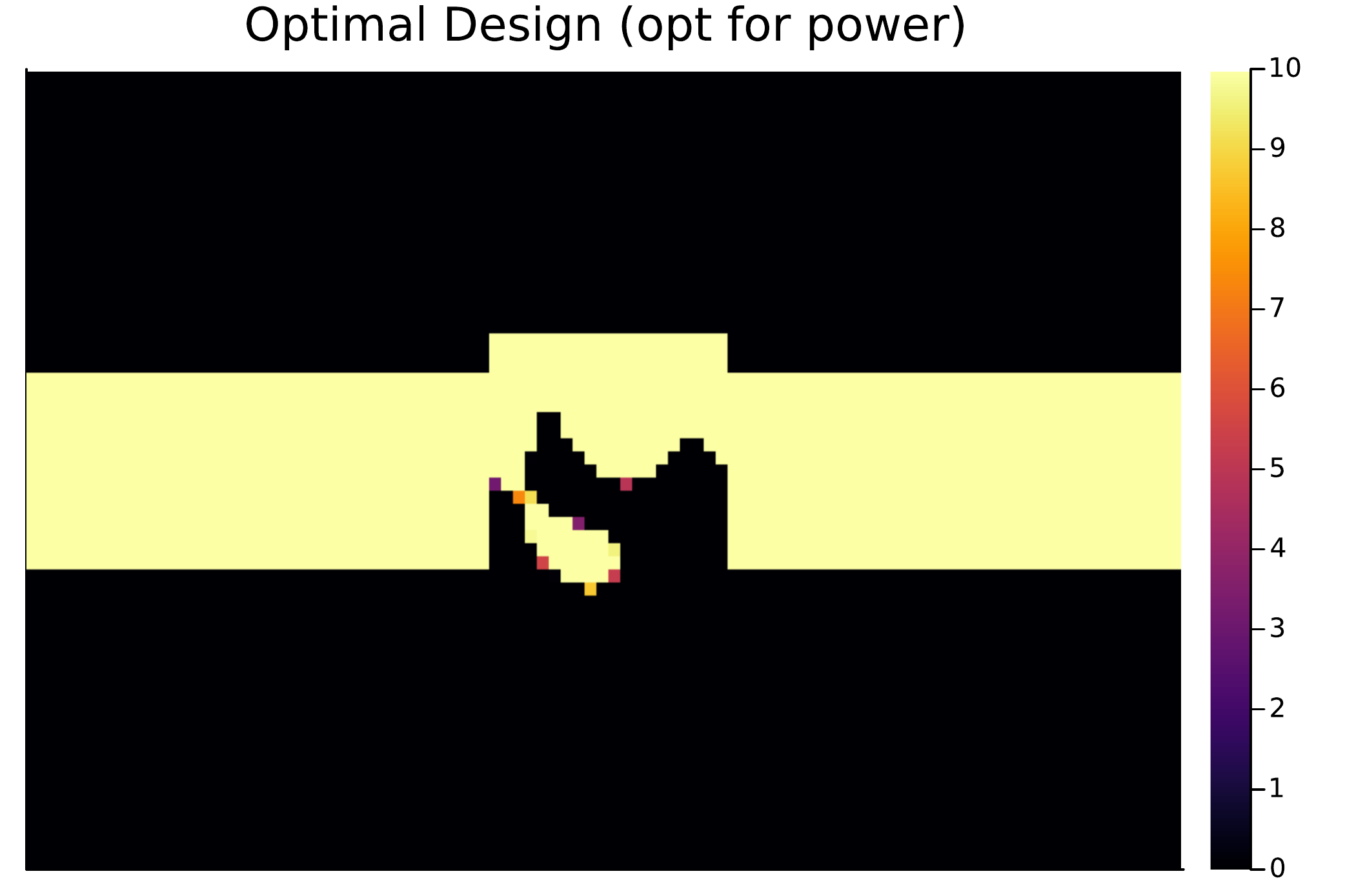}
    \includegraphics[width=0.45\linewidth]{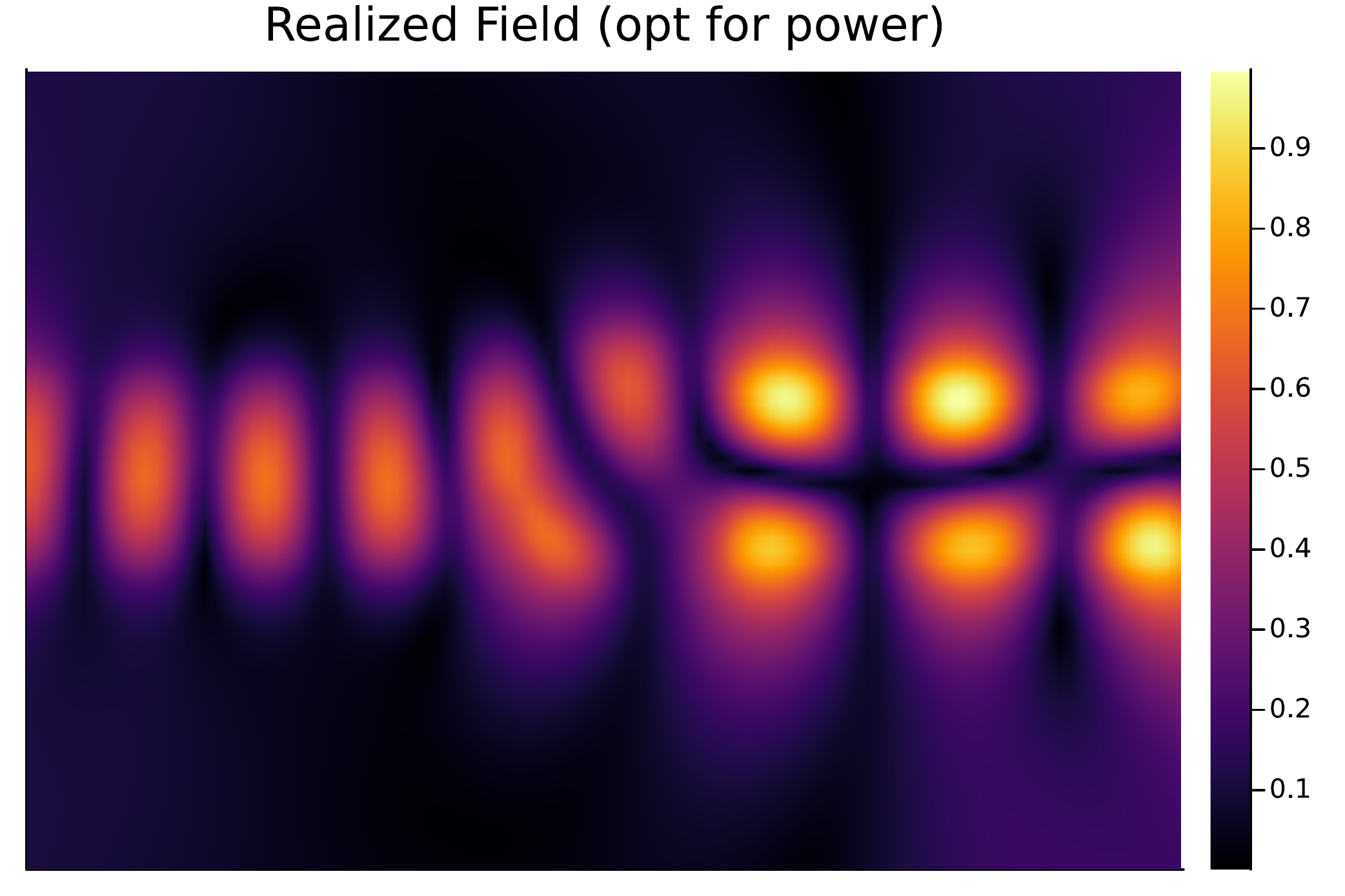}
    \caption{The design (left) is optimized for mode power. The corresponding field 
    (right) has greater power at the output compared to that of the 
    purity-optimized design, but it sacrifices some amount of purity.}
	\label{fig:power-design-field}
\end{figure}

\begin{figure}[t]
    \centering
	\includegraphics[width=0.7\linewidth]{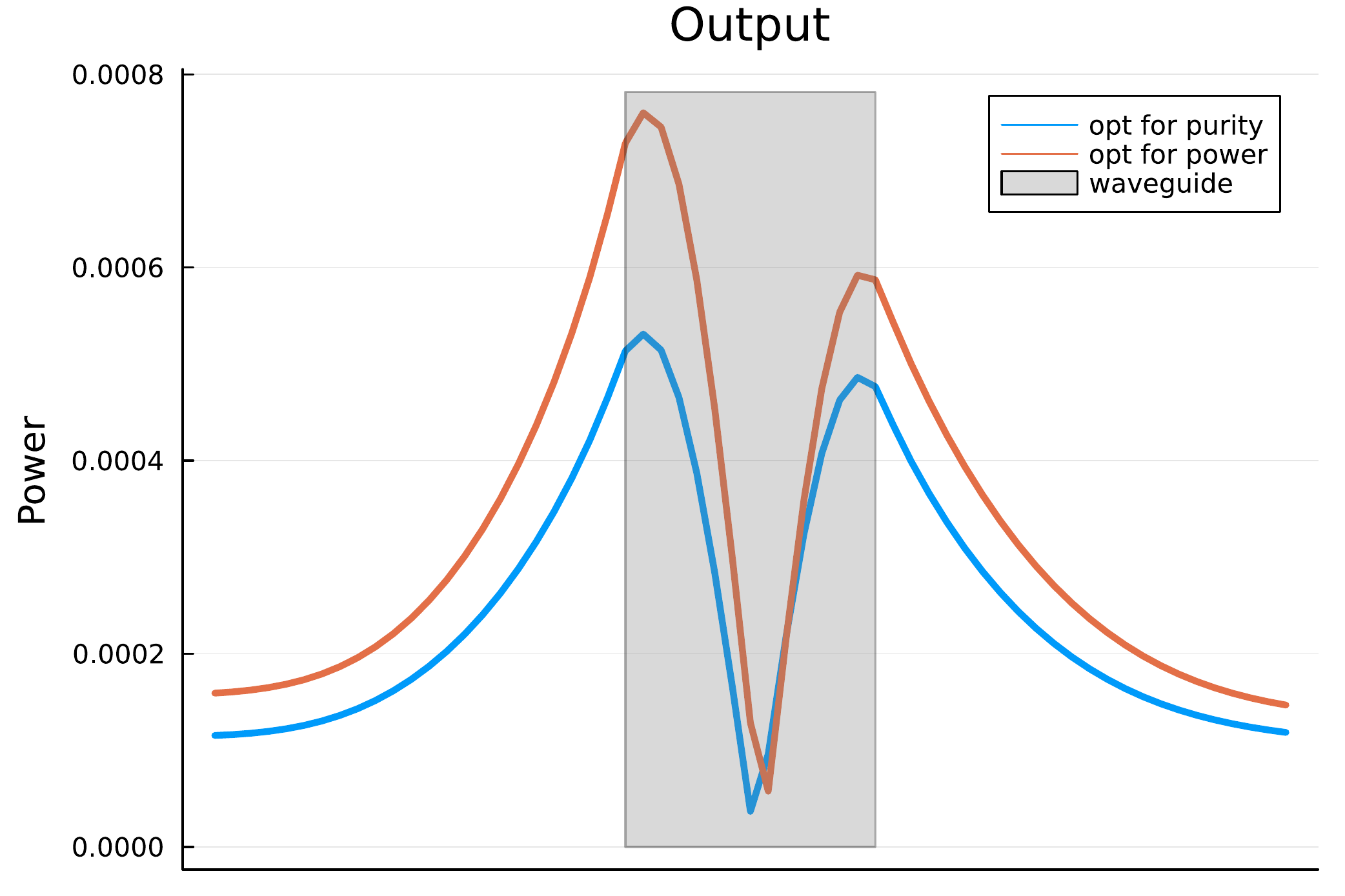}
    \caption{The design optimized for mode purity better matches the target
    mode waveform but has lower output power.}
	\label{fig:power-output}
\end{figure}

\paragraph{Results.} The resulting upper bound on the mode purity, that no design can
exceed, is .981. We also find an (approximately) optimal design with $\theta_i = \theta_{i+n}$
(\ie, with real permittivities). This design,
and its corresponding field, are shown in figure~\ref{fig:purity-deisgn-field}. The mode purity
this design achieves is .966, which is $(.981 - .966)/.981 \approx 1.5\%$ percent from the upper bound.
We note that this design, while very close to the optimal value for the mode purity,
is not very good in a practical sense: most of the power in the input waveguide is actually
scattered out to space. In general, we find that simply optimizing for the numerator,
as is usually done in practice, yields designs that are relatively efficient and have reasonable
mode purity. In this case, simply maximizing the numerator of the objective
results in a design that achieves a mode purity of .933, with an output power that is approximately 76\% greater.
(This design, and its corresponding field, is shown in figure~\ref{fig:power-design-field}.) 
This difference is highlighted in figure~\ref{fig:power-output}.

\section{Conclusion and future work}
In this paper, we have presented a simple method to compute
bounds on a number of efficiency metrics for physical design problems,
by solving a semidefinite program.
In particular, we focused on the common case where the efficiency metric
can be written as a ratio of two quadratics, which includes
metrics such as the focusing efficiency and the mode conversion efficiency.
We present a small example, but note that, while larger numerical examples
are possible, the resulting semidefinite programs are large; 
computing bounds on designs of larger sizes in reasonable time will likely 
require more sophisticated solvers (or larger computers). While the designs
shown here are also somewhat reasonable, they are still very far from the 
three dimensional designs that are useful in practice. Future work would
focus on creating faster solvers that can exploit the special structure of these problems,
along with simple interfaces that are user-friendly and can be used to easily
set up and solve these bounds.

\section*{Acknowledgements}
Theo Diamandis is supported by the Department of Defense (DoD) through the
National Defense Science \& Engineering Graduate (NDSEG) Fellowship Program.
The authors would also like to thank the anonymous reviewers for their comments
and suggestions, many of which we have incorporated in this text.
\bibliographystyle{ieeetr}
\bibliography{citations.bib}

\clearpage

\appendix
\section{Problem set up}\label{app:set-up}
The package uses an integral equation approximation to the Helmholtz equation
as the physics equation. We describe how the package solves this problem at a
high level in what follows.

\paragraph{Helmholtz's equation.} In this case, the initial physics equation is:
\[
\nabla^2 \phi(x) + k^2(1+\kappa(x))\phi(x) = f(x),
\]
for $x \in \Omega$. Here, $\Omega \subseteq \reals^2$ is a compact domain,
while $\phi: \Omega \to \complex$ is the (complex) amplitude of the field,
while $\kappa: \Omega \to [0, \kappa_\mathrm{max}]$ is the contrast, $k \in
\reals_+$ is the wavenumber, and $f: \Omega \to \complex$ is the excitation. We
will show that this can be approximated in the following form:
\[
z + G\diag(\theta)z = Gb,
\]
where $\theta_i = \kappa(x_i)/\kappa_\mathrm{max}$, $b_i = f(x_i)$, and  $z_i
\approx \phi(x_i)$ for some (chosen) points $x_i \in \Omega$. This is a common
method for computing approximate solutions to Helmholtz's equation
(\cf,~\cite[\S2.5]{petersonComputationalMethodsElectromagnetics1998}), but we
present it here for completeness.

\paragraph{Green's function.} Whenever $\kappa = 0$, \ie, when $\phi$ satisfies,
\[
\nabla^2 \phi(x) + k^2\phi(x) = f(x)
\]
there is a simple solution to the problem by a linear operator $\mathcal{G}$, such that
\[
\phi = \mathcal{G}f,
\]
where $\mathcal{G}$ is known as the \emph{Green's function} of the original equation:
\begin{equation}\label{eq:int}
(\mathcal{G}f)(x) = -\frac{\ii}{4}\int_\Omega H_0(k\|x-y\|)f(y)\,dy.
\end{equation}
Here, $H_0: \reals_+ \to \complex$ is the Hankel function of order zero of the first kind
(see, \eg,~\cite{boasMathematicalMethodsPhysical2006}).
Using this fact, we can then rewrite the original equation in terms of $\mathcal{G}$:
\[
\phi(x) + (\mathcal{G}(\kappa\phi))(x) = (\mathcal{G}f)(x),
\]
where $(\kappa \phi)(x) = \kappa(x)\phi(x)$ denotes the pointwise multiplication of the 
functions $\kappa$ and $\phi$.

\paragraph{Approximation.} We can then approximate the previous expression by taking a 
discretization. We assume that $\{x_1, \dots, x_n\} \subseteq \Omega$ denotes
a regularly-spaced grid with grid spacing $h > 0$. In this case, we will approximate the
equation in the following way:
\[
z + G\diag(\theta) z = Gb,
\]
where $z \in \complex^n$ is an approximation of the field amplitude $\phi$,
$G \in \complex^{n\times n}$ is the Green's operator, $b \in \complex^n$ is
the excitation, and $\theta \in [0, 1]^n$ are the permittivities
along the points of the grid. We can then make the following correspondences:
\[
\theta_i = \kappa(x_i)/\kappa_\mathrm{max}, \quad b_i = f(x_i), \quad i=1, \dots, n.
\]
while
\[
G_{ij} = -\delta (\ii/4) h H_0(k\|x_i - x_j\|), \quad i, j = 1, \dots, n, \quad i\ne j.
\]
This corresponds to approximating the integral~\eqref{eq:int} with a Riemann
sum on all of the off-diagonal terms. Because $H_0(0)$ is undefined, we will
approximate the diagonal terms of $G$ with the following integral:
\[
G_{ii} = -\frac{\ii\delta}{\pi}\int_0^{2\pi}\int_{0}^{h/2} tH_0(kt)\,dt\,d\rho = -4\ii \int_{0}^{h/2} tH_0(kt)\,dt = 
\frac{4}{\pi k^2} - \frac{\ii h}{k}H_1\left(\frac{kh}{2}\right),
\]
where $H_1$ is the Hankel function of order 1 of the first kind, while
$\delta = \kappa_\mathrm{max}$ is known as the maximum material contrast.
We can interpret this integral as integrating $H_0$ over a circle of radius $h/2$ and linearly
interpolating the resulting value to a square of side $h$ by scaling the result by
$h^2/(\pi (h/2)^2) = 4/\pi$.

With these definitions (and some additional regularity conditions on $f$, $\kappa$, and $\Omega$ 
which almost universally hold in practice) we then have that $z_i \approx \phi(x_i)$. In other
words, the solution to the discretized problem is approximately equal to the true solution
at the grid points $x_i$.

\section{Performance tricks}
In this section, we outline some additional tricks and tools which the overall
computation time of both the bounds and the heuristics when using this formulation.
Most of these ideas are implemented in whole or in part by the \texttt{WaveOperators.jl}
library, but we describe them here at a high level.

\paragraph{Removing zero-contrast points.} In many important practical cases,
we usually prefer to write the physics equation
\[
z + G\diag(\theta) z = Gb
\]
and constrain several entries of $\theta$ to be equal to zero (\ie, these entries imply
that there is no material present at position $x_i$ in the grid). 
In this case, it is
possible to separate $z$ into the components which have nonzero 
contrast $z_0$ and
positive contrast, $z_+$. We assume that the entries are in order such
that $z = (z_0, z_+)$ and $\theta = (0, \theta_+)$.
This means we can separate the physics equation into its
individual components
\[
G\diag(\theta)z =
\begin{bmatrix}
G_{00} & G_{0+}\\
G_{+0} & G_{++}    
\end{bmatrix}
\begin{bmatrix}
    \diag(\theta_0)z_0\\
    \diag(\theta_+)z_+
\end{bmatrix} =
\begin{bmatrix}
G_{00} & G_{0+}\\
G_{+0} & G_{++}    
\end{bmatrix}
\begin{bmatrix}
    0\\
    \diag(\theta_+)z_+
\end{bmatrix},
\]
where the diagonal matrices are square. Written out, after cancellations, we get
\[
\begin{aligned}
    z_0 + G_{0+}\diag(\theta_+)z_+ &= (Gb)_0\\
    z_+ + G_{++}\diag(\theta_+)z_+ &= (Gb)_+.
\end{aligned}
\]
Note that the first equation can be written as
\[
z_0 = (Gb)_0 - G_{0+}\diag(\theta_+)z_+ ,
\]
so no inverses need to be computed and all field
values $z_0$ can be easily written in terms of
the variables $\diag(\theta_+)$ and $z_+$ only,
while $z_+$ does not depend on the values of $z_0$.

\paragraph{Schur complement.} In some other special cases, it is also
easier to specify parameters $\theta_+$ which might be nonzero,
but are fixed ahead of time. For convenience, we will write
$\theta_+ = (\theta_c, \theta_f)$, where $\theta_c$ are the 
nonzero parameters that are constrained, while $\theta_f$ are the
free parameters, and similarly for $z_+ = (z_c, z_f)$. In this case, we can similarly
separate the physics equation into its individual components:
\[
G_{++}\diag(\theta_+)z_+ = \begin{bmatrix}
    G_{cc} & G_{cf}\\
    G_{fc} & G_{ff}
\end{bmatrix}\begin{bmatrix}
    \diag(\theta_c)z_c\\
    \diag(\theta_f)z_f
\end{bmatrix}.
\]
This results in the following physics equations over the points with nonzero contrast:
\[
\begin{aligned}
    z_c + G_{cc}\diag(\theta_c)z_c + G_{cf}\diag(\theta_f)z_f &= (Gb)_c\\
    z_f + G_{fc}\diag(\theta_c)z_c + G_{ff}\diag(\theta_f)z_f &= (Gb)_f.
\end{aligned}
\]
We can then eliminate the variable $z_c$ form this equation to receive a linear
equation that depends only on the product of the free parameters and the
points corresponding to the free field, $\diag(\theta_f)z_f$. To do this,
we solve for $z_c$ in the first equation:
\[
    z_c = (I+G_{cc}\diag(\theta_c))^{-1}((Gb)_c - G_{cf}\diag(\theta_f)z_f),
\]
and plug it into the second to get
\[
z_f + (G_{ff} + \bar G_{cc})\diag(\theta_f)z_f = (Gb)_f + \bar b_c,
\]
with
\[
\bar G_{cc} = -G_{fc}(I+G_{cc}\diag(\theta_c))^{-1}G_{cf}, \qquad \bar b_c = -G_{fc}(I+G_{cc}\diag(\theta_c))^{-1}(Gb)_c,
\]
which is easily seen to be of the form of~\eqref{eq:greens_formalism}. Because the
SDP size scales quadratically on the number of field variables, and
SDPs themselves usually have a large runtime, this is often
a useful procedure as it only requires computing the matrix $\bar G_{cc}$ once at the
beginning of the problem. This then reduces the total number of variables in the SDP at
the expense of computing a single matrix factorization at the beginning of the procedure.

\section{Pareto Frontier}
We plot the Pareto frontier for the problem of optimizing mode purity and output
power in a mode converter, considered in this paper's numerical experiments.
\begin{figure}
    \centering
    \includegraphics[width=0.7\linewidth]{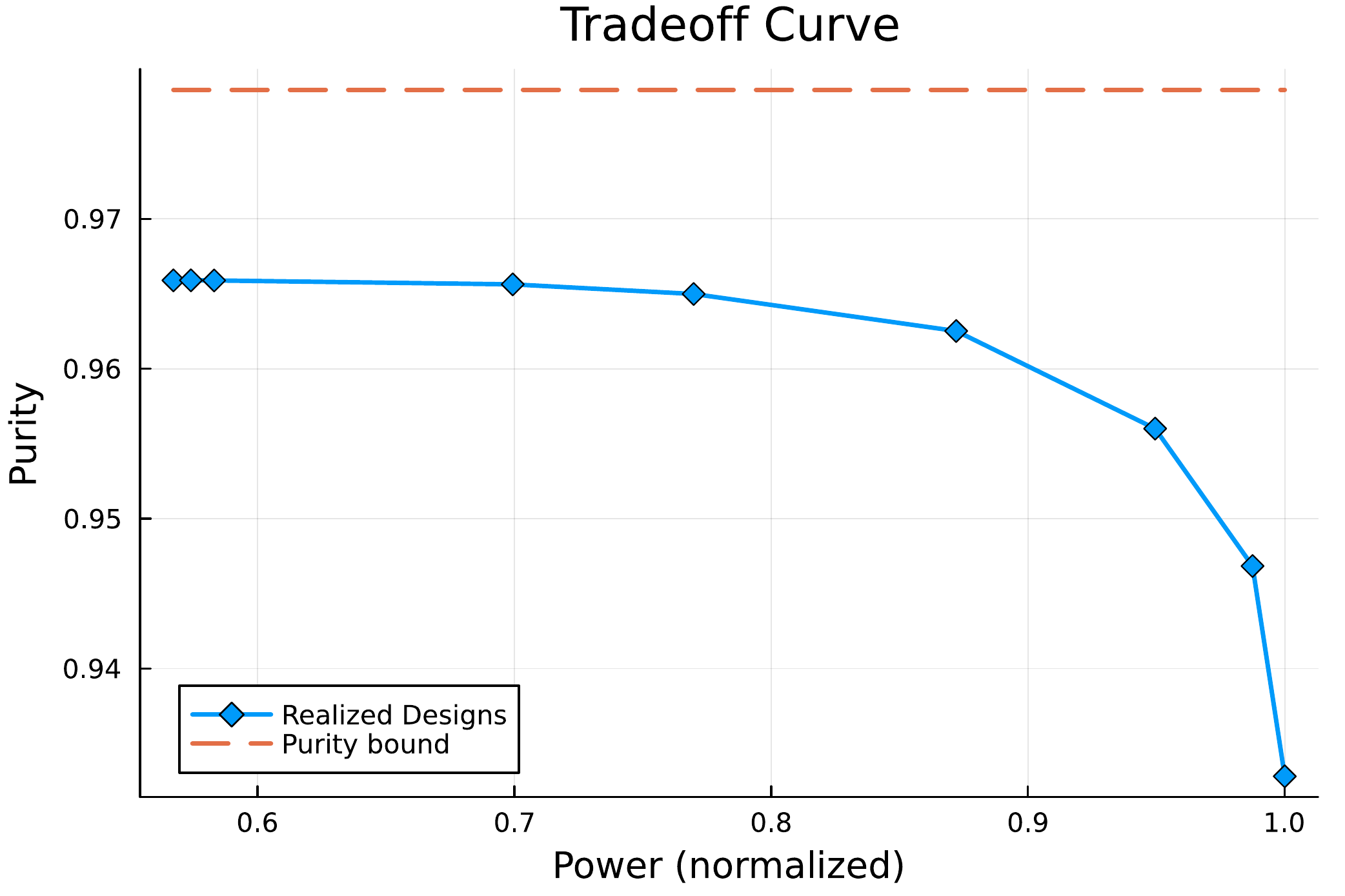}
    \caption{The Pareto frontier for the problem of optimizing mode purity and output
    power. Each point on the plot indicates a realized design. The SDP bound on 
    mode purity (with no constraint on power) is also shown.}
	\label{fig:pareto}
\end{figure}

\end{document}